\documentclass{ifacconf}
\pdfoutput=1

\usepackage{amsmath}
\usepackage{amssymb}
\usepackage{enumerate}
\usepackage{overpic}
\usepackage{xcolor}

\usepackage{graphicx}      
\usepackage{natbib}        
\usepackage{siunitx}
\usepackage{makecell}
\usepackage{multirow}
\usepackage[hyphens]{url}

\allowdisplaybreaks[4]
\def\fuel{\text{fuel}}
\def\gt{\text{gt}}
\def\el{\text{em}}
\def\bat{\text{b}}
\def\ch{\text{c}}
\def\drv{\text{drv}}
\def\CDCS{\text{CDCS}}

\newcommand{\vvec}[1]{\overrightarrow{#1}}

\begin{document}
\begin{frontmatter}

\title{Optimal energy management for hybrid electric aircraft} 


\author[first]{Martin Doff-Sotta}
\quad
\author{Mark Cannon\,$\mbox{}^\ast$}
\quad
\author{Marko Bacic\,$\mbox{}^{\ast,1}$}

\address[first]{Department of Engineering Science, University of Oxford, UK\\ Email: \{martin.doff-sotta, mark.cannon, marko.bacic\}@eng.ox.ac.uk}

\thanks{On part-time secondment from Rolls-Royce PLC.} 

\begin{abstract}                
A convex formulation is proposed for optimal energy management in aircraft with hybrid propulsion systems consisting of gas turbine and electric motor components. 
By combining a point-mass aircraft dynamical model with models of electrical and mechanical powertrain losses, the fuel consumed over a planned future flight path is minimised subject to constraints on the battery, electric motor and gas turbine. The resulting optimisation problem is used to define a predictive energy management control law that takes into account the variation in aircraft mass during flight.
A simulation study based on a representative 100-seat aircraft with a prototype parallel hybrid electric propulsion system is used to investigate the properties of the controller. 
We show that an optimisation-based control strategy can provide significant fuel savings over heuristic energy management strategies in this context.
\end{abstract}

\begin{keyword}
Energy Management, Nonlinear Model Predictive Control, Convex Programming.
\end{keyword}

\end{frontmatter}

\section{Introduction}
Aviation currently contributes around 2\% of worldwide human-made CO2 emissions but the demand for air travel and transport is growing at a significant rate. The aviation industry is committed to realising this growth sustainably with a drastic reduction of CO2 emissions by 2050. One avenue identified to contribute to the required CO2 reduction is through hybridisation of aircraft propulsion systems. This refers to enabling technologies for boundary layer ingesting aircraft \citep{hall2017} as well as  rotary/tilt wing aircraft configurations in the Urban Air Mobility markets \citep{nasa_urban}. Hybrid electric architectures require real-time dynamic power management in order to minimise CO2 output. 

This paper addresses an optimal energy management problem for a hybrid electric aircraft with a propulsion system consisting of a gas turbine and a battery-powered electric motor in a parallel configuration. 
Although we consider here a battery as a secondary energy source, the approach is equally applicable to other primary and secondary energy sources such as hydrogen powered reciprocating engines, fuel-cells and super-capacitors.

Any optimisation methodology for primary power management must satisfy the basic requirements of determinism, convergence in finite time and verifiability. 
We propose a solution based on model predictive control (MPC) employing convex optimisation. Predicted performance, expressed in terms of the fuel consumption over a given future flight path, is optimised subject to constraints on power flow and stored energy, and subject to the nonlinear aircraft dynamics, which include nonlinear losses in powertrain components. The proposed convex formulation of the optimisation problem is made suitable for a real time nonlinear MPC implementation by introducing several key simplifying assumptions on the characteristics of powertrain components. Specifically, the gas turbine and the electric motor are modelled via sets of convex quasi-static power maps, battery losses are modelled using a time-invariant equivalent circuit, and the available data on the future flight path is assumed sufficient to determine powertrain shaft speeds across the prediction horizon. 

Supervisory control methodologies for energy management have been proposed in the context of hybrid electric ground vehicles~\citep[e.g.][]{sciarretta07}.  
Several approaches have been proposed for this problem, including methods based on indirect optimal control~\citep{kim11,onori15},  Dynamic Programming~\citep{lin03} and MPC~\citep{koot05,east19ieeetcst}. 
Optimal control of hybrid propulsion systems in aircraft is a new application area that poses a number of distinct challenges, perhaps the most significant of which are complex nonlinear flight dynamics and the effects of the time-varying aircraft weight due to the burning of fuel during flight. On the other hand, the future power demand is likely to be more reliably predictable in aircraft than in cars since a pre-planned flight path is generally available for aircraft, whereas the driving cycle is subject to greater uncertainty in route and traffic conditions~\citep{dicairano14,josevski17}.
The contribution of this paper is to demonstrate that the optimal energy management problem for hybrid aircraft can be posed as a convex optimisation problem.
To the authors' knowledge, this is the first attempt to address an important new application area of energy management.

The paper is organised as follows. Section~\ref{sec:modelling} derives a continuous time hybrid electric aircraft model. This model is the basis of the discrete-time model and the MPC strategy that are proposed in Section~\ref{sec:dt_model}. Section~\ref{sec:convex} shows that the minimisation of fuel consumption can be expressed as a convex problem. Section~\ref{sec:simulation_results} describes simulation results and
conclusions are drawn in Section~\ref{sec:conclusion}.



\section{Modelling}\label{sec:modelling}

\begin{figure}[h!]
\centerline{\begin{overpic}[scale=0.644,trim=0mm 0mm 0mm -2mm,clip]{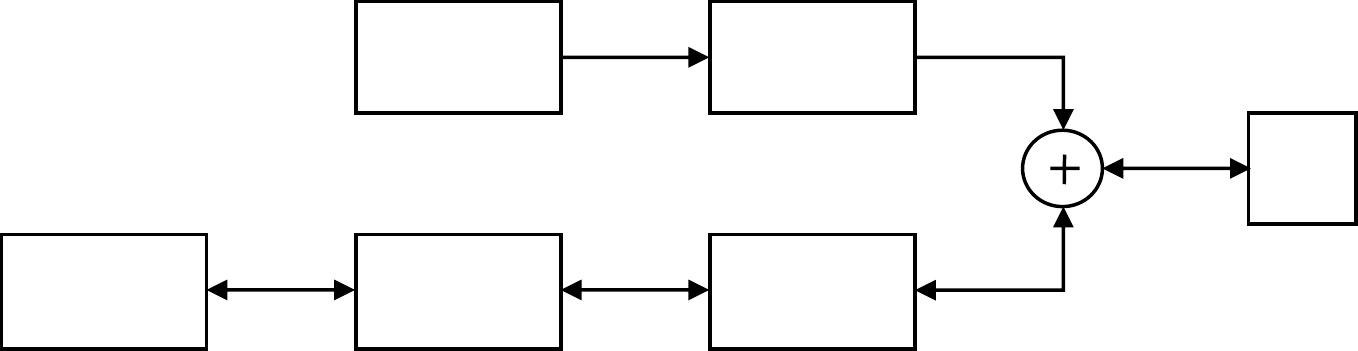}
\put(1.9,3.7){\footnotesize{Battery}}
\put(28.6,5){\footnotesize{Electric}}
\put(31.5,1.7){\footnotesize{Bus}}
\put(54.5,5){\footnotesize{Motor/}}
\put(52.4,1.7){\footnotesize{Generator}}
\put(30.8,20.6){\footnotesize{Fuel}}
\put(56.7,22.3){\footnotesize{Gas}}
\put(54,19.2){\footnotesize{Turbine}}
\put(93.4,12.5){\footnotesize{Fan}}
\put(18.5,6.5){$P_{\bat}$}
\put(44.5,6.5){$P_{\ch}$}
\put(70,6.5){$P_{\el}$}
\put(42.5,23.8){$P_{\fuel}$}
\put(70,23.8){$P_{\gt}$}
\put(83.1,15.5){$P_{\drv}$}
\end{overpic}}
    \caption{Hybrid-electric propulsion architecture.}
    \label{fig:propulsion}
\end{figure}
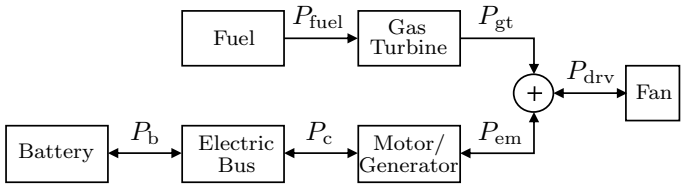

We assume a parallel hybrid electric aircraft propulsion system in which a gas turbine producing power $P_{\gt}(t)$ is combined with an electric motor with power output $P_{\el}(t)$ (Fig.~\ref{fig:propulsion}). The net power output of the propulsion system, $P_{\drv}(t)$, is produced by combining these two power sources via the relation $P_{\drv}(t) = P_{\gt}(t) + P_{\el}(t)$ (assuming $100\%$ efficiency in drivetrain components). When the drive power is negative, which may occur for example while the aircraft is descending, it is assumed that the same powertrain could be used to generate electrical energy (i.e.\ it is capable of operating in a ``windmilling'' mode) in order to recharge the battery. In practice, a variable-pitch fan would be required, which increases complexity. The gas turbine and electric motor shaft rotation speeds are $\omega_{\gt}(t)$ and $\omega_{\el}(t)$ respectively. 

The electric motor is powered by a battery with state of charge (SOC) $E(t)$ and rate of change of energy $P_{\bat}$, and the state of charge dynamics are given by 
\[
\dot{E}=-P_{\bat}.  
\]
The battery is modelled as an equivalent circuit with internal resistance $R$ and open-circuit voltage $U$, so that
\begin{align*}
P_{\bat}&=g\bigl(P_{\el}\left(t\right), \omega_{\el}\left(t\right)\bigr) \\
&=\frac{U^2}{2R}\left(1-\sqrt{1-\frac{4R}{U^2}h(P_{\el}(t), \omega_{\el}(t)) }\right). 
\end{align*}
where $U$ and $R$ are assumed constant~\citep{east2018}.
The map relating the mechanical power output of the electric motor $P_{\el}$ to electrical input power $P_{\ch}$ is $P_{\ch} = h(P_{\el})$. 
We assume that, for fixed $\omega_{\el}$, $h(P_{\el},\omega_{\el})$ is non-decreasing and differentiable with respect to $P_{\el}$ and $h(\cdot)$ is determined empirically from electric motor loss map data as
\[
h(P_{\el}, \omega_{\el}) = \kappa_{2}(\omega_{\el})P_{\el}^2  + \kappa_{1}(\omega_{\el})P_{\el}  + \kappa_{0}(\omega_{\el})
\]
for some functions $\kappa_2(\cdot)$, $\kappa_1(\cdot)$, $\kappa_0(\cdot)$.
with $\kappa_{2}(\omega_{\el}) \geq 0$ and $\kappa_{1}(\omega_{\el}) > 0$ for all $\omega_{\el}$ in the operating range. 

\begin{figure}[htb]
\centerline{\begin{overpic}[scale=0.35,trim=-1mm 0mm 0mm 8mm, clip]{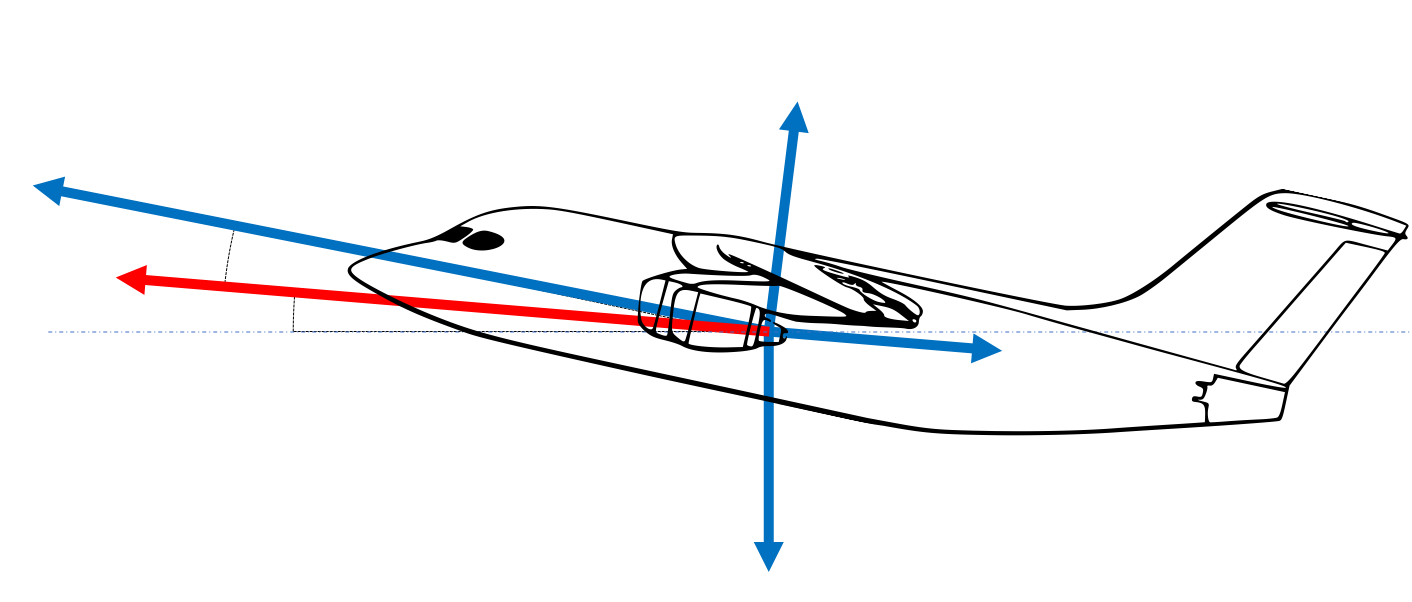}
\put(-1.5,28){$\vvec{T}$}
\put(72,14){$\vvec{D}$}
\put(56.5,2){$\vvec{W}$}
\put(54.5,37){$\vvec{L}$}
\put(4,22){$\vvec{v}$}
\put(17.5,23.3){$\alpha$}
\put(17.5,19){$\gamma$}
\end{overpic}}
\vspace{-2mm}
\caption{Aircraft forces and motion.}
\label{fig:aircraft}
\end{figure}
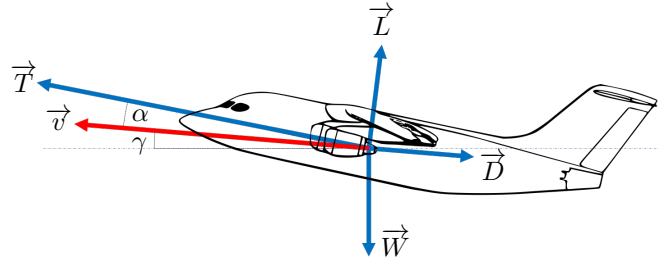

The aircraft motion is constrained by its dynamic equations. Assuming a point-mass model~\citep{stevens2015aircraft} and referring to Figure \ref{fig:aircraft}, the equilibrium of forces yields
\begin{equation}
m\frac{\mathrm{d}}{\mathrm{d}t}( \vvec{v})  =\vvec{T} + \vvec{L} + \vvec{D} + \vvec{W}
\label{eq:T_t}
\end{equation}
where $ \vvec{v}$ is the velocity vector, $m$ the instantaneous mass of the aircraft, $\vvec{T}$ the vector of thrust, $\vvec{L}$ and $\vvec{D}$ are the lift and drag vectors and $\vvec{W}$ is the aircraft weight. 

Using the polar coordinates parametrisation ($v$,$\gamma$), the drive power is given as follows
\[
P_{\drv}=\vvec{T} \cdot \vvec{v} = m\frac{\mathrm{d}}{\mathrm{d}t}(\frac{1}{2} v^2)  + \frac{1}{2}C_D\rho S v^3 +  m g \sin{(\gamma)} v 
\]
where $v$ is the magnitude of the velocity vector, $S$ is the wing area, $\rho$ the density of air, $g$ the acceleration due to gravity, $\gamma$ the flight path angle, $C_D=C_D(\alpha)$ the drag coefficient, $\alpha$ the angle of attack.
Similarly, projecting equation \eqref{eq:T_t} along the lift vector $\vvec{L}$) yields
\[
m v \frac{\mathrm{d}}{\mathrm{d}t}( \gamma) + mg \cos{(\gamma)}= T \sin{(\alpha)} + \frac{1}{2}C_L\rho S v^2 ,
\]
%
where $C_L=C_L(\alpha)$ is the lift coefficient. The contribution of the thrust in the vertical direction being very small, the term $T \sin{(\alpha)}$ can be neglected (it can be checked a posteriori that $\alpha$ is small).

The rate of change of the aircraft mass is given by
\[
 \dot{m}  = \dot{m}_{\fuel} = -f(P_{\gt}(t), \omega_{\gt}(t))
\]
where $\dot{m}_{\fuel}$ is the rate of fuel consumption and $f(P_{\gt},\omega_{\gt})$ is assumed to be convex, differentiable and non-decreasing with respect to $P_{\gt}$ for fixed $\omega_{\gt}$. We assume that $f(\cdot)$ can be determined empirically from fuel map data in the form 
\[
f(P_{\gt},\omega_{\gt}) = \beta_{2}(\omega_{\gt})P_{\gt}^2 + \beta_{1}(\omega_{\gt})P_{\gt}  + \beta_{0}(\omega_{\gt}) , 
\]
with $\beta_{2}(\omega_{\gt})\!\geq\! 0$, $\beta_{1}(\omega_{\gt}) \!>\! 0$ in the operating range of $\omega_{\gt}$. 

The problem at hand is to find the real-time optimal power split between the gas turbine and electric motor that minimises
\begin{equation}\label{eq:contin-time_objective}
    J = \int_{0}^{T}{f(P_{\gt}(t), \omega_{\gt}(t))}\mathrm{d}t
\end{equation}
while satisfying constraints on the battery SOC, limits on power flows throughout the powertrain, and producing sufficient power to follow a prescribed flight path. 




\section{Discrete-time optimal control}\label{sec:dt_model}

This section describes a discrete-time model that enables the optimisation of the power split between the electric motor and the gas turbine over a given future flight path to be formulated as a finite-dimensional optimisation problem.
For a fixed sampling interval $\delta$, we consider a predictive control strategy that minimises, online at each sampling instant, the predicted fuel consumption over the remaining flight path. 
The minimisation is performed subject to the discrete-time dynamics of the aircraft mass and the battery SOC. The problem is also subject to bounds on the stored energy in the battery (to prevent deep discharging or overcharging), as well as limits on power flows that correspond to physical and safety constraints.

The optimal solution to the fuel minimisation problem at the $k$th sampling instant is computed using estimates of the battery SOC $E(k\delta)$ and the aircraft mass $m(k\delta)$. The control law at time $k\delta$ is defined by the first time step of this optimal solution.
The notation $\{ x_{0} , x_{1} , \ldots x_{N-1}\}$ is used for the sequence of future values of a variable $x$ predicted at the $k$th discrete-time step, so that $x_i$ denotes the predicted value of $x\bigl((k+i)\delta\bigr)$. The horizon $N$ is chosen so that $N = \lceil T/\delta \rceil - k$, and hence $N$ shrinks as $k$ increases and $k\delta$ approaches $T$.

The discrete-time approximation of the objective (\ref{eq:contin-time_objective}) is  
\begin{equation}
    J = \sum_{i=0}^{N-1}{f(P_{\gt,i},\omega_{\gt, i})} \, \delta
\label{eq:obj}
\end{equation}
with, for $i=0,\dots,N-1$,
\begin{align}
\mbox{}\hspace{-3mm} f(P_{\gt, i}, \omega_{\gt, i}) &= \beta_{2}(\omega_{\gt})P_{\gt,i}^2 \!+\! \beta_{1}(\omega_{\gt, i})P_{\gt, i}  \!+\! \beta_{0}(\omega_{\gt, i})
\label{eq:m_k0}
\\
m_{i+1} &= m_{i} - f(P_{\gt,i}, \omega_{\gt,i}) \, \delta \label{eq:m_k}
\end{align}
where the forward Euler approximation has been used for derivatives. Using the same approach, the discrete-time battery model is
\begin{align}
E_{i+1} &= E_{i} - g(P_{\el,i}, \omega_{\el,i}) \, \delta
\label{eq:E_k}
\\
 P_{\bat,i} &= g (P_{\el,i}, \omega_{\el, i}) 
\nonumber 
\\
&= \frac{U^2}{2R}\biggl[1-\sqrt{1-\frac{4R}{U^2}h(P_{\el, i}, \omega_{\el, i}) }\biggr]
\label{eq:P_b}
\end{align}
for $i=0,\ldots,N-1$, where
\begin{multline}\label{eq:h_k}
h(P_{\el, i},\omega_{\el, i}) =\\
\kappa_{2}(\omega_{\el, i})P_{\el, i}^2 
+ \kappa_{1}(\omega_{\el, i})P_{\el, i} + \kappa_{0}(\omega_{\el, i})
\end{multline}
and the aircraft dynamics are given in discrete time by
\begin{align}
& 
m_{i} v_{i} \Delta_{i} \gamma = - m_{i} g \cos{(\gamma_{i})} + \tfrac{1}{2}C_L(\alpha_{i})\rho S v_{i}^2
\label{eq:vertical_k}
\\
& 
P_{\drv,i} = \begin{aligned}[t] &\tfrac{1}{2} m_{i} \Delta_{i} (v^2)  +  m_{i} g \sin{(\gamma_{i})}v_{i}\\
&+\tfrac{1}{2}C_D(\alpha_{i})\rho S v_{i}^3,  
\end{aligned}
\label{eq:drive}
\\
& 
P_{\drv,i} =P_{\gt,i} + P_{\el,i}, 
\label{eq:P_link}
\end{align}
for $i=0,\ldots,N-1$, where 
\[
\Delta_{i} (v^2) = (v^2_{i+1}   -   v^2_{i})/\delta , \quad
\Delta_{i} \gamma = (\gamma_{i+1}   -   \gamma_{i})/\delta.
\]
%
The problem to be solved at each time step $k$ is therefore: 
\begin{align}
& \min_{\substack{P_{\gt},\,P_{\el},\,P_{\drv}\\m,\,E,\,\omega_{\gt},\,\omega_{\el}}}
& & \sum^{N-1}_{i=0} f(P_{\gt,i}, \omega_{\gt,i})
\label{eq:min} \\
& \qquad \text{ s.t.} & &  P_{\drv,i} = P_{\gt,i} + P_{\el,i} 
\nonumber \\
& & &   \begin{aligned} P_{\drv,i} &= \tfrac{1}{2} m_{i} \Delta_{i} v^2   +  m_{i} g \sin{(\gamma_{i})}v_{i} \\ &\quad +\tfrac{1}{2}C_D(\alpha_{i})\rho S v_{i}^3 \end{aligned}
\nonumber\\
& & &   m_{i} v_{i} \Delta_{i} \gamma =  - m_{i} g \cos{(\gamma_{i})} + \tfrac{1}{2}C_L(\alpha_{i})\rho S v_{i}^2 
\nonumber\\
& & &   m_{i+1}  =m_{i} - f(P_{\gt,i}, \omega_{\gt,i}) \,\delta
\nonumber\\
& & &   E_{i+1} = E_{i} - g(P_{\el,i},\omega_{\el,i}) \, \delta
\nonumber\\
& & &   m_{0} = m(k\delta) 
\nonumber\\
& & &   E_{0} = E(k\delta) 
\nonumber\\
& & &  \underline{E} \leq E_{i} \leq  \overline{E}
\nonumber\\
& & &  \underline{P}_{\gt} \leq P_{\gt,i} \leq  \overline{P}_{\gt}
\nonumber\\
& & &  \underline{\omega}_{\gt} \leq \omega_{\gt,i} \leq  \overline{\omega}_{\gt}
\nonumber\\
& & &  {\underline{P}_{\el}} \leq P_{\el,i} \leq  \overline{P}_{\el}
\nonumber\\
& & &  {\underline{\omega}_{\el}} \leq \omega_{\el,i} \leq  \overline{\omega}_{\el}
\nonumber
\end{align}
where the constraints are imposed for $i=0,\ldots,N-1$.
Here ${(\overline{E},\underline{E})}$ are the bounds on  SOC that are required for normal battery operation, ${(\overline{P}_{\gt},\underline{P}_{\gt})}$ and ${(\overline{P}_{\el}, \underline{P}_{\el})}$ are the bounds on gas turbine power and electric motor power respectively, and ${ (\overline{\omega}_{\el}, \underline{\omega}_{\el})}$ and ${(\underline{\omega}_{\gt},\overline{\omega}_{\gt})}$ are the bounds on the gas turbine and electric motor shaft rotation speeds.
%




\section{Convex formulation}\label{sec:convex}

The optimisation problem in (\ref{eq:min}) is nonconvex, which makes a real-time implementation of an MPC algorithm that relies on its solution computationally intractable. In this section a convex formulation is proposed that is suitable for an online solution.
We assume that the aircraft speed $v_i$ and flight path angle $\gamma_i$ are chosen externally by a suitable guidance algorithm for $i=0,\ldots,N-1$.

A convex formulation of the drive power is derived by expressing the drag and lift coefficients, $C_D$ and $C_L$, as functions of the angle of attack $\alpha$ and combining the equations that constrain the aircraft motion in the forward and vertical directions. Over a restricted domain and for given Reynolds and Mach numbers, the drag and lift coefficients can be expressed respectively as a quadratic non-decreasing function and a linear non-decreasing function \citep{abbott1945}:
\begin{alignat}{2}
C_D(\alpha_{i}) &= a_2 \alpha_{i}^2 + a_1 \alpha_{i} + a_0, &\qquad  a_2 &> 0
\label{eq:C_D}
\\
C_L(\alpha_{i}) &= b_1 \alpha_{i} + b_0, &\qquad  b_1 &> 0
\label{eq:C_L}
\end{alignat}
for $\underline{\alpha} \leq \alpha_i \leq \overline{\alpha}$.

Combining \eqref{eq:vertical_k}, \eqref{eq:drive}, \eqref{eq:C_D} and \eqref{eq:C_L}, the angle of attack can be eliminated from the expression for drive power, so that $P_{\drv,i}$ can be expressed as a quadratic function of the aircraft mass, $m_{i}$, as follows
\begin{equation}\label{eq:P_drv}
    P_{\drv,i} = \eta_{2,i} m_{i}^2 + \eta_{1,i} m_{i} + \eta_{0,i} , 
\end{equation}
where
\begin{align*}
&    \eta_{2,i} = \frac{ 2 a_2 (v_{i} \Delta_{i} \gamma + g \cos{(\gamma_{i})}  )^2 }{ b_1^2 \rho S v_{i}} ,
\\
&    \eta_{1,i} = \tfrac{1}{2}\Delta_{i} v^2 + g \sin{(\gamma_{i})} v_{i}  \\
&\qquad\quad - \frac{2 a_2 b_0  (v_{i} \Delta_{i} \gamma
    + g \cos{(\gamma_{i})} ) v_{i}}{b_1^2} 
\\
&\qquad\quad    + \frac{a_1}{ b_1 }  ( v_{i} \Delta_{i} \gamma + g \cos{(\gamma_{i})}  ) v_{i} ,
\\
&    \eta_{0,i} = \tfrac{1}{2} \rho S v_{i}^3 \Bigl(  \frac{a_2 b_0^2}{b_1^2} - \frac{a_1 b_0}{b_1} + a_0  \Bigr), 
\end{align*}
Here the flight path angles $\gamma_i$ and speeds $v_i$ are assumed to be fixed and are not optimisation variables. Since $\eta_{2,i} > 0$ for all $i$, the drive power is a convex function of  $m_{i}$. Note that there is no guarantee that satisfying equation \eqref{eq:P_drv} enforces equations \eqref{eq:vertical_k} and \eqref{eq:drive} individually.  In practice,  assuming that we have full control over the eliminated variable $\alpha$ (via the elevator and fans), both individual dynamical equations can be satisfied a posteriori. The inequality constraint on $\alpha$ also has to be checked a posteriori. 

For the given parallel hybrid configuration, we assume for simplicity that the electric motor and gas turbine share a common shaft rotation speed which is equal to the speed of rotation of the fan, i.e.\ $\omega_{\gt,i} = \omega_{\el,i}$ for all $i$. If the shaft speed is known at each discrete-time step of the prediction horizon, then the coefficients in \eqref{eq:m_k0} and \eqref{eq:h_k} can be estimated from a set of polynomial approximations of $h(\cdot)$ and $f(\cdot)$ at a pre-determined set of speeds. This allows $h(P_{\el,i},\omega_{\el,i})$ and $f(P_{\gt,i},\omega_{\gt,i})$ in (\ref{eq:P_b}) and (\ref{eq:min}) to be replaced by time-varying functions of the gas turbine power and electric motor power alone:
\begin{align}
h(P_{\el,i},\omega_{\el,i}) &= h_i(P_{\el,i}) 
= \kappa_{2,i}P_{\el,i}^2  + \kappa_{1,i}P_{\el,i}  + \kappa_{0, i},
\label{eq:h_P}
\\
f(P_{\gt,i},\omega_{\gt,i}) &= f_i(P_{\gt,i}) 
= \beta_{2,i} P_{\gt,i}^2 + \beta_{1,i} P_{\gt,i}  + \beta_{0, i} 
\label{eq:f_P}
%
\end{align}
with $\kappa_{2,i} \geq 0$, $\kappa_{1,i} >0$ and $\beta_{2,i} \geq 0$, $\beta_{1,i} > 0$ for all $i$.

In order to estimate the shaft speed $\omega_{\gt,i} = \omega_{\el,i}$, and hence determine the coefficients $\kappa_{2,i}$, $\kappa_{1,i}$, $\kappa_{0,i}$, $\beta_{2,i}$, $\beta_{1,i}$, $\beta_{0,i}$ in (\ref{eq:h_P}) and (\ref{eq:f_P}), we use a pre-computed look-up table relating the drive power to rotational speed of the fan, for a given altitude,  Mach number and air conditions (temperature and specific heat at constant pressure). This enables the shaft speed to be determined as a function of the fan power output at each discrete-time step along the flight path. Although $P_{\drv,i}$ depends on the aircraft mass $m_i$, which is itself an optimisation variable, a prior estimate of the required power output can be obtained from (\ref{eq:P_drv}) assuming a constant mass $m_i=m_0$ for all $i$. The simulation results described in Section~\ref{sec:simulation_results} show that this assumption has a negligible effect on solution accuracy. 

We define $g_i(\cdot)$ in terms of $h_i(\cdot)$ as 
\[
g_i(P_{\el,i})
=\frac{U^2}{2R}\biggl[ 1-\sqrt{1-\frac{4R}{U^2}h_i(P_{\el, i}) }\biggr] .
\]
Then $g_i(\cdot)$ is necessarily a convex, non-decreasing, one-to-one function if the lower bound on $P_{\el,i}$ is redefined as
%
\[
\underline{P}_{\el,i} :=\max \Bigl\{ -\overline{P}_{\el}, -\frac{\kappa_{1,i}}{2 \kappa_{2,i}}\Bigr\},
\]
since this bound ensures that $h_i(\cdot)$ is a one-to-one non-decreasing convex function of  $P_{\el,i}$. 

The dynamic constraints  \eqref{eq:m_k}, \eqref{eq:E_k} and the power balance (\ref{eq:P_link}) can be expressed using (\ref{eq:P_drv}), (\ref{eq:h_P}) and (\ref{eq:f_P}) as 
\begin{align}
&    m_{i+1} = m_{i} - f_i(P_{\gt,i}) \, \delta
   \label{eq:m_k_P}
   \\
&    E_{i+1} = E_{i} - g_i (P_{\el,i})  \, \delta 
    \label{eq:E_k_P}
\\
& P_{\gt,i} + P_{\el,i} = \eta_{2,i} m_i^2 + \eta_{1,i} m_i + \eta_{0,i} .
\label{eq:drv_k_P}
\end{align}
These constraints are nonconvex due to their quadratic dependence on
the optimisation variables $P_{\gt,i}$, $P_{\el,i}$ and $m_i$.
To convexify these constraints, we first eliminate $P_{\el,i}$ from (\ref{eq:E_k_P}) and (\ref{eq:drv_k_P}) using $P_{\bat,i} = g_i(P_{\el,i})$ and $P_{\el,i} = g_i^{-1}(P_{\bat,i})$. Then \eqref{eq:E_k_P} becomes linear,
\[
    E_{i+1} = E_{i} - P_{\bat,i} \delta  .
\]
Moreover, under the assumptions on $g_i(\cdot)$ (convex, non-decreasing and one-to-one), the inverse mapping $g_i^{-1}(\cdot)$
is a concave, increasing function \citep[e.g.][]{east2018}. Note that $g^{-1}(\cdot)$ is given explicitly as
\[
    g_i^{-1}(P_{\bat,i})=-\frac{\kappa_{1,i}}{2\kappa_{2,i}} + \biggl[ -\frac{R P^2_{\bat,i}}{\kappa_{2,i}U^2} + \frac{P_{\bat,i} - \kappa_{0,i}}{\kappa_{2,i}} +\frac{\kappa^2_{1,i}}{4\kappa^2_{2,i}}\biggr]^{\frac{1}{2}} 
\]
if $\kappa_{2,i} > 0$,
%
and by
\[
    g_i^{-1}(P_{\bat,i})=-\frac{1}{\kappa_{1,i}} \Bigl( \frac{R}{U^2}P^2_{\bat,i} - P_{\bat,i} + \kappa_{0,i} \Bigr)
\]
at any time steps $i$ such that $\kappa_{2,i}=0$.
Therefore, by relaxing the equality constraints in (\ref{eq:m_k_P}) and (\ref{eq:drv_k_P}) to inequalities, a pair of convex constraints: 
\begin{align}
& m_{i+1} \leq m_{i} - f_i(P_{\gt,i}) \, \delta
\label{eq:m_ineq_k_P}
\\
& P_{\gt,i} \geq \eta_{2,i} m_i^2 + \eta_{1,i} m_i + \eta_{0,i} - g_i^{-1}(P_{\bat,i}) 
\label{eq:drv_ineq_k_P}
\end{align}
is obtained since $g_i^{-1}(\cdot)$ is concave and $f_i(\cdot)$ is convex. 

With these modifications, and noting that the objective in (\ref{eq:min}) is equivalent to minimising $m_0-m_ N$, the optimisation to be solved to determine the optimal power profile at the $k$th time step can be expressed as the convex problem:
\begin{align}
& \min_{\substack{P_{\gt},P_{\bat}\,P_{\drv}\\m,\,E,\,\omega_{\gt},\,\omega_{\el}}}
& & m_{0} - m_N
\label{eq:min2} \\
& \qquad \text{ s.t.} & & P_{\gt,i} \geq \eta_{2,i} m_{i}^2 + \eta_{1,i} m_{i} + \eta_{0,i} - g_i^{-1}(P_{\bat,i})
\nonumber \\
& & &   m_{i+1}  \leq m_{i} - f_i (P_{\gt,i}) \, \delta
\nonumber \\
& & &   E_{i+1} = E_{i} - P_{\bat,i}\, \delta 
\nonumber \\
& & &   m_{0} = m(k\delta)
\nonumber \\
& & &   E_{0} = E(k\delta) 
\nonumber \\
& & &  \underline{E} \leq E_{i} \leq  \overline{E}
\nonumber \\
& & &  \underline{P}_{\gt} \leq P_{\gt,i} \leq  \overline{P}_{\gt}
\nonumber \\
& & &  \underline{P}_{\bat,i} \leq P_{\bat,i} \leq  \overline{P}_{\bat,i}
\nonumber
\end{align}
where $\underline{P}_{\bat,i}=g_i(\underline{P}_{\el,i})$, $\overline{P}_{\bat,i}=g_i(\overline{P}_{\el,i})$, and the constraints are imposed for $i=0,\ldots,N-1$. The form of the objective in (\ref{eq:min2}) ensures that any feasible solution that does not satisfy the constraints in (\ref{eq:m_ineq_k_P}) and (\ref{eq:drv_ineq_k_P}) with equality is suboptimal. 
%
Thus the solutions of (\ref{eq:min2}) and (\ref{eq:min}) are necessarily equal if (\ref{eq:min}) is feasible.







\section{Numerical results}
\label{sec:simulation_results}

This section uses the optimisation problem \eqref{eq:min2} to construct an energy management case study involving a representative hybrid-electric passenger aircraft.
Solutions of \eqref{eq:min2} were computed using the general purpose convex optimisation solver CVX~\citep{cvx}. Since the minimisation in (\ref{eq:min2}) is convex, convergence of the solver to a global optimum is ensured. 

\subsection{Simulation scenario}\label{sec:sim_setup}

The parameters of the model used in simulations are collected in Table \ref{tab:param}. These are based on publicly available data for the BAe 146 aircraft.
The propulsion system is assumed to consist of four gas turbines and electric motors, each with the hybrid-parallel configuration shown in Fig.~\ref{fig:propulsion}. 

\begin{table}[h!]
\begin{tabular}{llll}
\hline
\textbf{Parameter} & \textbf{Symbol} & \textbf{Value} & \textbf{Units} \\ \hline
Mass (MTOW)     &    $m$    &   $42000$    &    \si{kg}   \\ \hline
Gravity acceleration     &    $g$    &   $9.81$    &    \si{m.s^{-2}}   \\ \hline
Wing area &   $S$     &    $77.3$   &    \si{m^2}   \\ \hline
Density of air &   $\rho$     &    $1.225$   &    \si{kg.m^{-3}} \\ \hline
\multirow{2}*{Lift coefficients} &   $b_0$   &    $0.43$   &    \si{-}  \\
&   $b_1$     &    $0.11$   &    \si{deg^{-1}}  \\ \hline
 \multirow{3}*{Drag coefficients} &   $a_0$   &   $0.029$   &    \si{-}  \\
&   $a_1$   &   $0.004$   &    \si{deg^{-1}} \\
&   $a_2$   &   $5.3\mathrm{e}{-4}$   &    \si{deg^{-2}}\\ \hline
Angle of attack range &   $\left[\underline{\alpha}; \overline{\alpha}\right]$     &    $\left[-3.9; 10\right]$   &    \si{deg}   \\ \hline
Fuel mass &   $m_{\fuel}$     &    $8000$   &    \si{kg}   \\ \hline
\multirow{2}*{Fuel map coefficients} &   $\beta_0$     &    $ 0.03$    &   \si{kg.s^{-1}} \\ 
&   $\beta_1$     &    $0.08$    &   $\si{kg.{MJ}^{-1}}$ \\  \hline
Battery SOC range &   $\left[\underline{E}; \overline{E}\right]$     &    $\left[221; 939\right]$   &    \si{MJ} \\ \hline
Gas turbine power range &   $\left[\underline{P}_{\gt}; \overline{P}_{\gt}\right]$     &    $\left[0; 5\right]$  &    \si{MW}   \\ \hline
Motor power range &   $\left[\underline{P}_{\el}; \overline{P}_{\el}\right]$     &    $\left[0; 2\right]$   &    \si{MW}   \\ \hline
$\#$ of arrangements &   $n$     &    $4$   &    \si{-} \\ \hline
Flight time &   $T$     &    $3600$   &    \si{s} \\ \hline
\end{tabular}
\vspace{1mm}\caption{Model parameters.}
\label{tab:param}
\vspace{-2mm}
\end{table}

\begin{figure}
\centerline{\includegraphics[width = .45\textwidth]{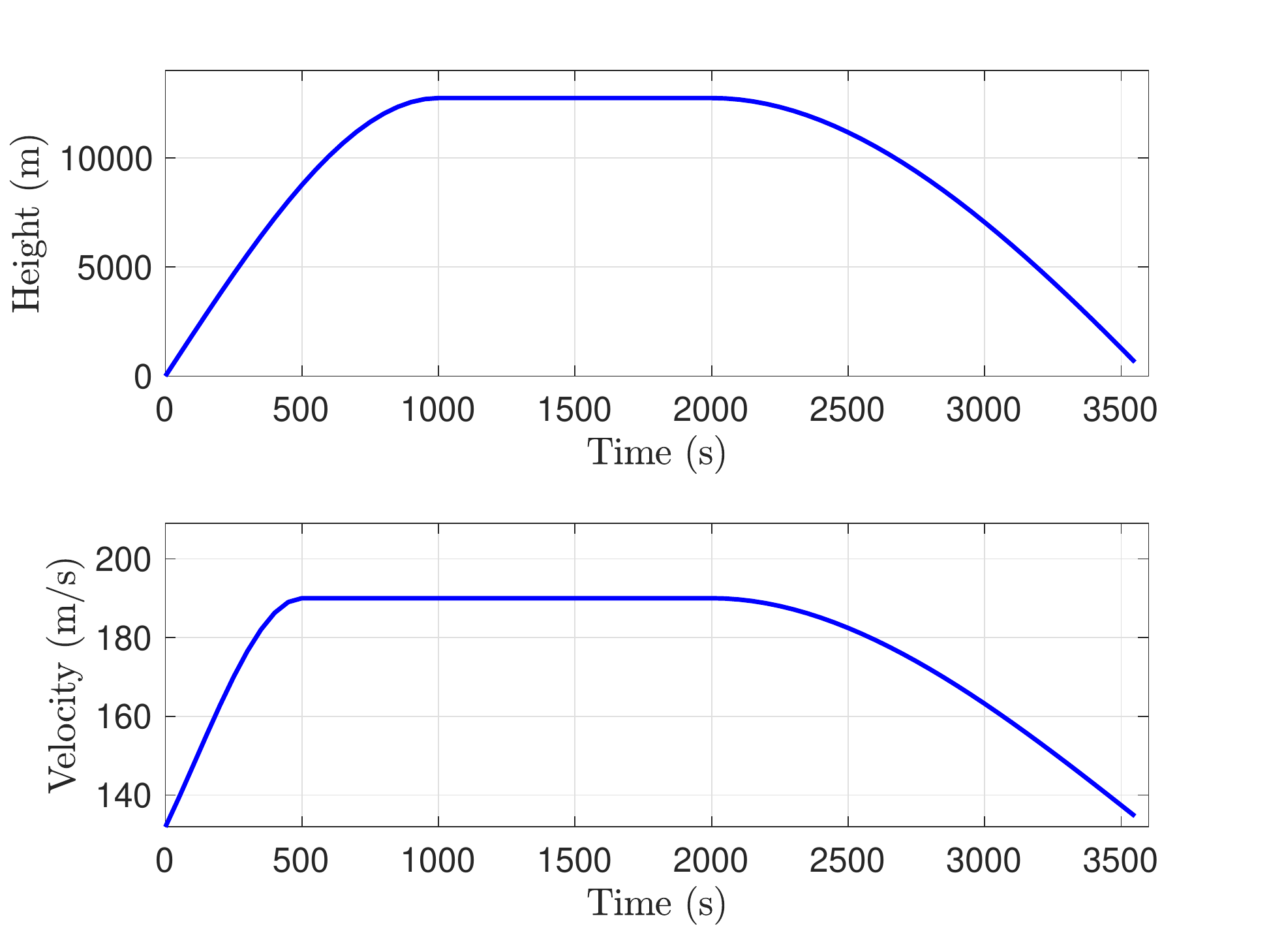}}
\caption{Height and velocity profiles for the mission.}
\label{fig:profile}
\end{figure}

For the purposes of this study it is assumed that velocity and height profiles are known a priori as a result of the fixed flight plan entered prior to take-off.  
We consider an exemplary 1-hour flight at a true airspeed (TAS) of $190$\,$\si{m/s}$ for a typical 100-seat passenger aircraft. The flight path (height and velocity profile) is shown in Figure \ref{fig:profile}. 

\begin{figure}
\includegraphics[width = .45\textwidth]{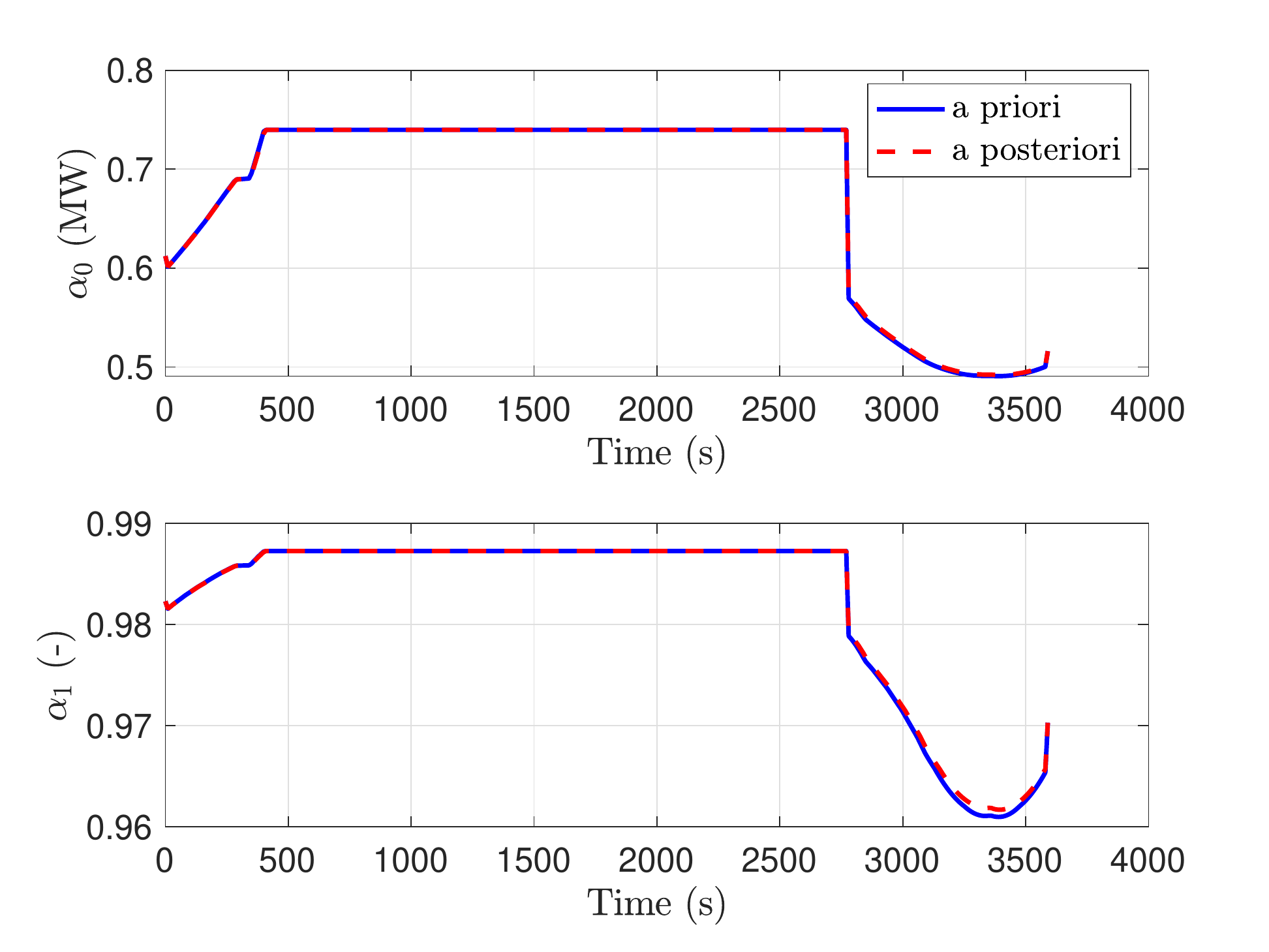}
\caption{Electric loss map coefficients computed with estimated drive power and actual drive power.}  
\label{fig:alpha}
\end{figure}

\begin{figure}
\includegraphics[width = .45\textwidth]{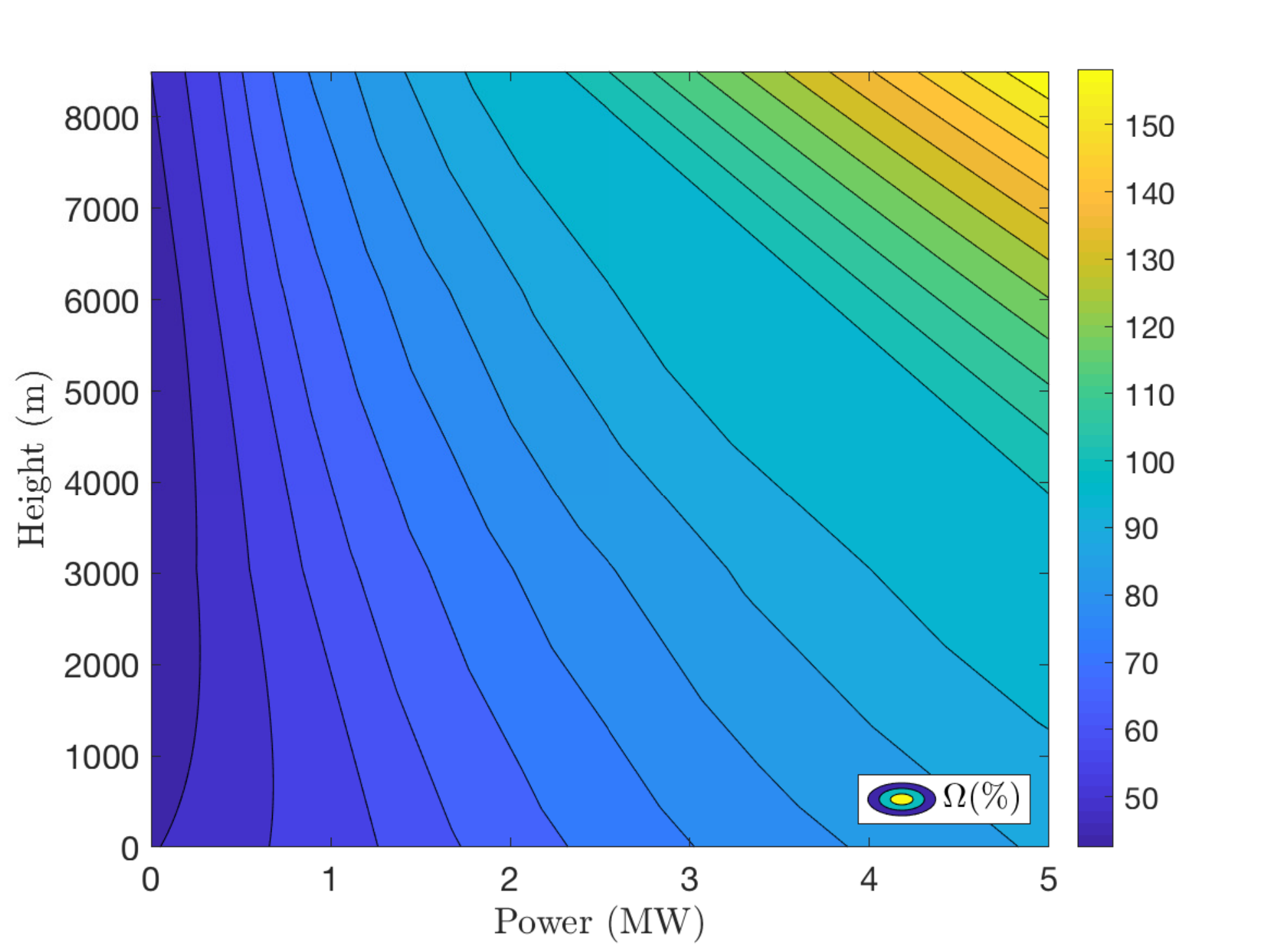}
\caption{Contour plot relating drive power, altitude and non-dimensional rotation speed ($\Omega$) for a Mach number of $0.55$. The shaft rotation speed is obtained from $\Omega$ as $\omega = \smash{\frac{156.7}{100}\frac{\pi}{30} \Omega \sqrt{T_{in}}}$ where $T_{in}=T_{0}(h) + v^2/2 c_p$ is the temperature at inlet of the fan, $c_p=1000$ $\si{J K^{-1} kg^{-1}}$ is the specific heat of air at constant pressure and $T_0(h)$ is the temperature of air at altitude $h$.}
\label{fig:fanmap}
\end{figure}

The electric loss map coefficients $\kappa_{i,2},\kappa_{i,1},\kappa_{i,0}$ can be estimated in two steps from these profiles. First, the shaft rotation speed, $\omega_i$ ($=\omega_{\gt,i} = \omega_{\el,i}$), is interpolated from a precomputed look-up table relating measured shaft rotation speed, altitude and drive power at a given Mach number (Fig.~\ref{fig:fanmap}). Then, the coefficients are interpolated from a precomputed record of losses in the electric motor as a function of rotation speed. This procedure requires drive power $P_{\drv}$ to be approximated \textit{a priori}, e.g.~by assuming constant aircraft mass for the duration of the flight. This assumption is supported by Figure~\ref{fig:alpha}, which shows that the electric map coefficients are almost identical for the estimated drive power profile and for the actual drive power profile computed retrospectively. We also find that the $\kappa_{2, i}$ coefficients are negligible for all $i$.
 
The gas turbine fuel map used in this study is approximately linear ($\beta_{2,i} \approx 0 \, \forall i$) and furthermore the fuel consumption does not depend significantly on shaft rotation speed.
Therefore the fuel map coefficients are given in Table \ref{tab:param} as constants (i.e.\ $\beta_{1,i} = \beta_1, \, \beta_{0,i} = \beta_0 \, \forall i$).


\subsection{Results}

\begin{figure}
\includegraphics[width = .45\textwidth]{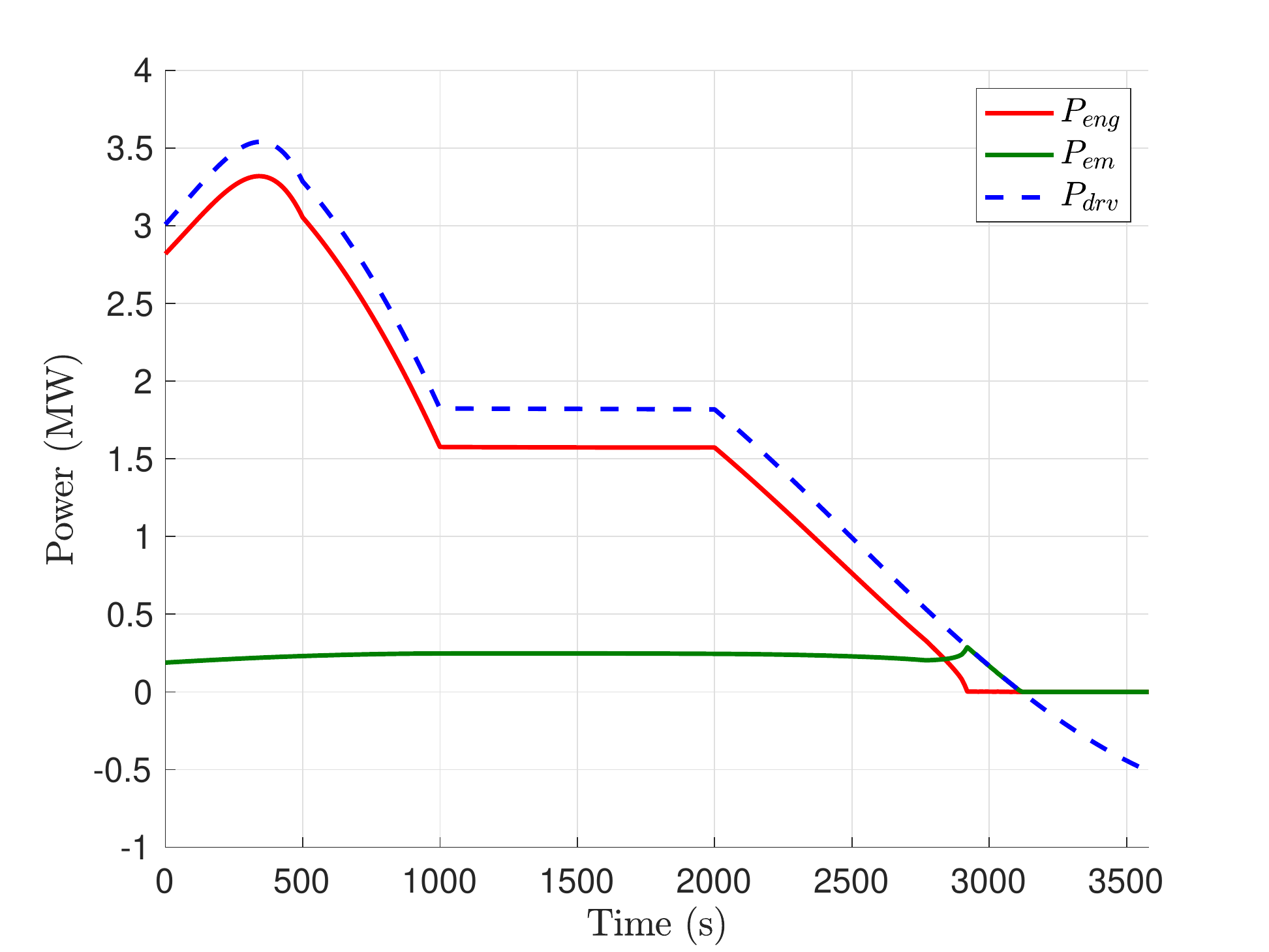}
\caption{MPC power split strategy obtained by solving \eqref{eq:min2} at each sampling instant with a shrinking horizon.}
\label{fig:Psplit}
\end{figure}

\begin{figure}
\includegraphics[width = .45\textwidth]{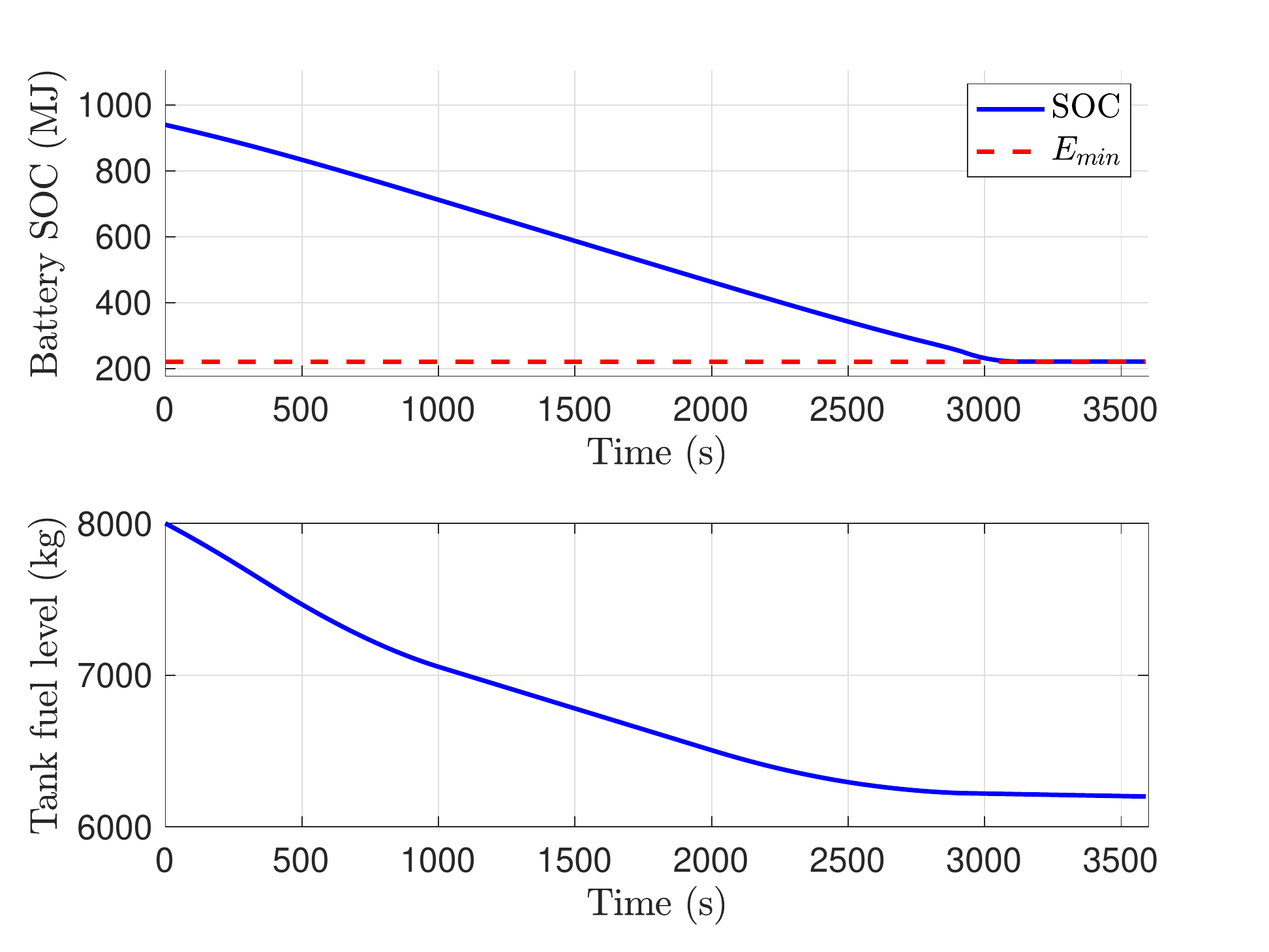}
\caption{Closed loop evolution of SOC and fuel mass.}
\label{fig:SOC}
\end{figure}

The mission is simulated with sampling interval $\delta=10$ $\si{s}$ over a one-hour shrinking horizon by solving the optimisation problem \eqref{eq:min2} at each time step and implementing the first element of the optimum power split sequence as an MPC law. 
The closed loop energy management control strategy is shown in Figure \ref{fig:Psplit}, which gives the power split for a single coupled gas turbine and electric motor.
%
Clearly the constraints on the gas turbine and electric motor power are respected.
The evolution of the battery SOC and fuel consumption are shown in Figure \ref{fig:SOC}. The upper plot illustrates that the constraints on SOC are respected and that the SOC reaches a minimum when the drive power becomes negative, as expected. The lower plot in Fig.~\ref{fig:SOC} shows that, as expected, the rate of fuel consumption is greater during the initial climb phase when the gas turbine power output is high. The fuel consumption recorded for this simulation is $F^{\ast}=1799$ $\si{kg}$. In comparison, a fully gas turbine-powered flight with the same initial total aircraft weight  has a fuel consumption of
$F_{\gt}=2034$ $\si{kg}$. We note however that this reduction is achieved at the expense of reduced available payload as a result of the weight of the electric components of the powertrain (battery storage and electric motors).

\begin{figure}
\includegraphics[width = .45\textwidth]{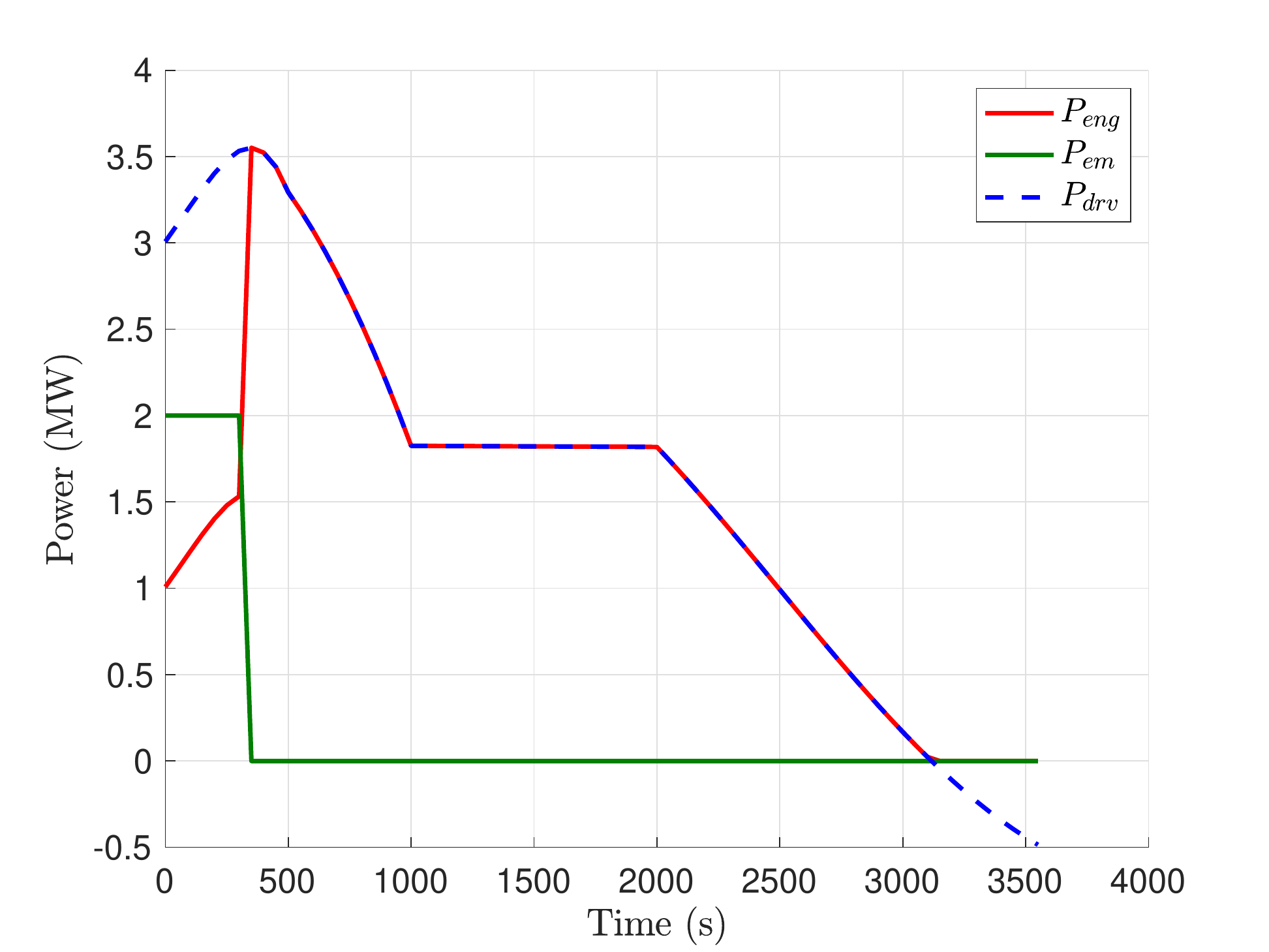}
\caption{Power split with a CDCS strategy.}
\label{fig:CDCS}
\end{figure}

In order to evaluate the optimality of the power split solution, we compare it with the strategy of supplementing the gas turbine with the maximum electric motor power ($\overline{P}_{\el}$) until the battery is fully depleted, then switching to a sustaining mode in which only the gas turbine operates. In hybrid vehicles this is known as a Charge-Depleting-Charge-Sustaining (CDCS) strategy~\citep{onori15}. Using this strategy the power split is as shown in Figure \ref{fig:CDCS} and the fuel consumption is $F_{\CDCS} = 1858 \si{kg}$. 

\begin{figure}
\includegraphics[width = .45\textwidth]{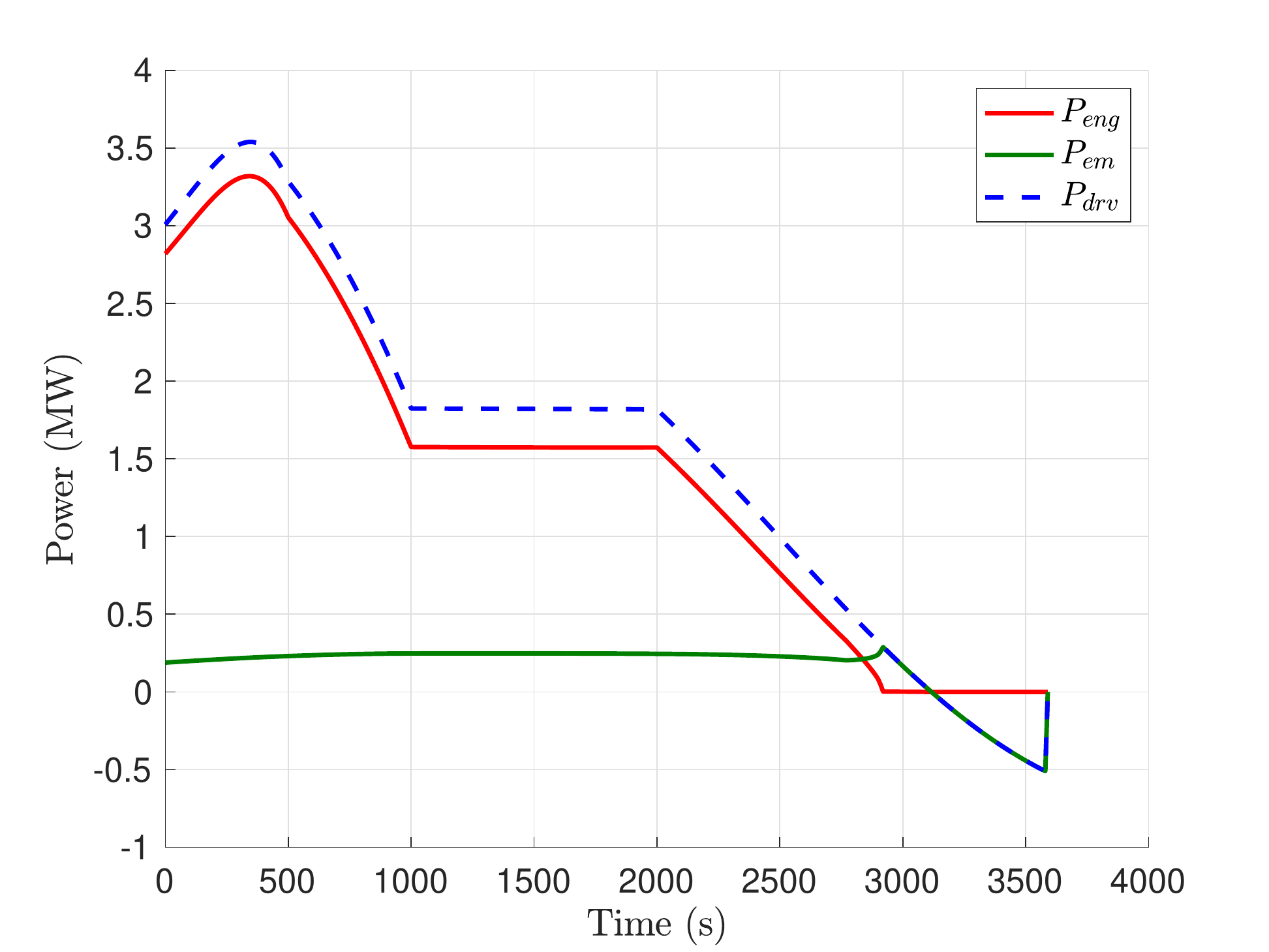}
\caption{Optimal power split with windmilling.}
\label{fig:Psplit4}
\end{figure}

\begin{figure}
\includegraphics[width = .45\textwidth]{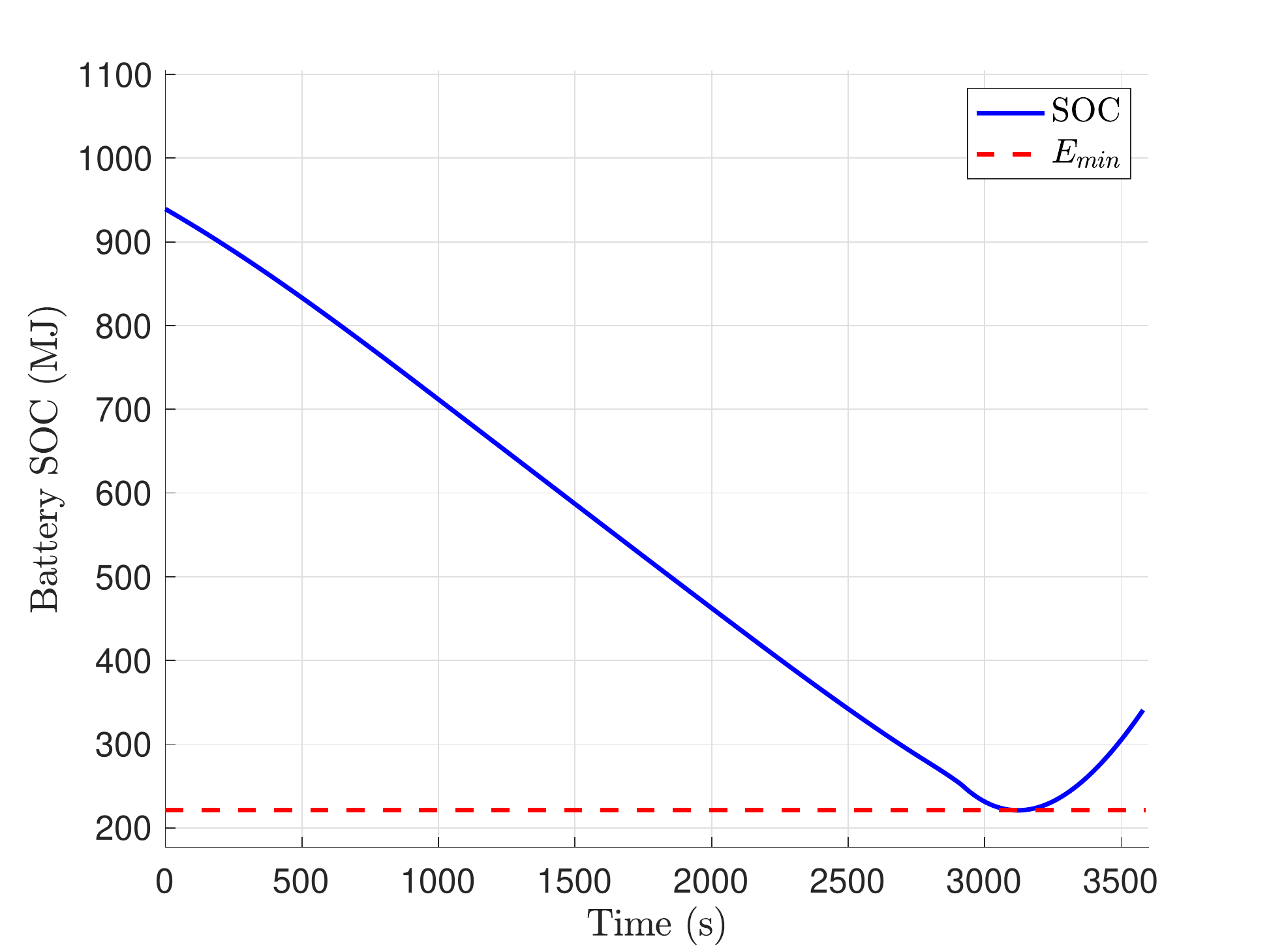}
\caption{Battery state of charge with windmilling.}
\label{fig:SOC4}
\end{figure}

To investigate the potential for windmilling (energy recovery when the net power demand is negative), the lower bound on electric motor power is set to $\underline{P}_{\el}=-2$ $\si{MW}$, to allow transmission of power from the fan to the battery with the electric motor acting as a generator. The optimisation problem (\ref{eq:min2}) is also modified
by introducing a terminal term in the objective function so as to maximise the SOC of the battery at the end of the flight: $J = m_0 - m_N - \lambda E_N.$ The coefficient $\lambda$ should be small to avoid adversely affecting the main objective of minimising fuel consumption. 
%
Replacing the objective for $\lambda = 0.1$ gives the results shown in Figures~\ref{fig:Psplit4} and \ref{fig:SOC4}.
The windmilling effect can be seen at the end of the flight and is characterised by negative electrical power and battery recharge.

\subsection{Discussion}

Intuitively, the optimal power split strategy might be expected to consume as much fuel as possible at the beginning of the flight so as to reduce the aircraft mass, and thus reduce the drive power needed during level flight and descent. However, the MPC strategy maintains an almost constant electric power over the whole flight (Fig.~\ref{fig:Psplit}). This is explained by the relatively short flight duration and the characteristics of the aircraft model, as a result of which the change in total mass of the aircraft is relatively small (less than $5\%$). Despite this, the MPC strategy achieves a non-negligible reduction in fuel consumption ($3.2\%$) over the CDCS strategy. 

\begin{figure}
\includegraphics[width = .45\textwidth]{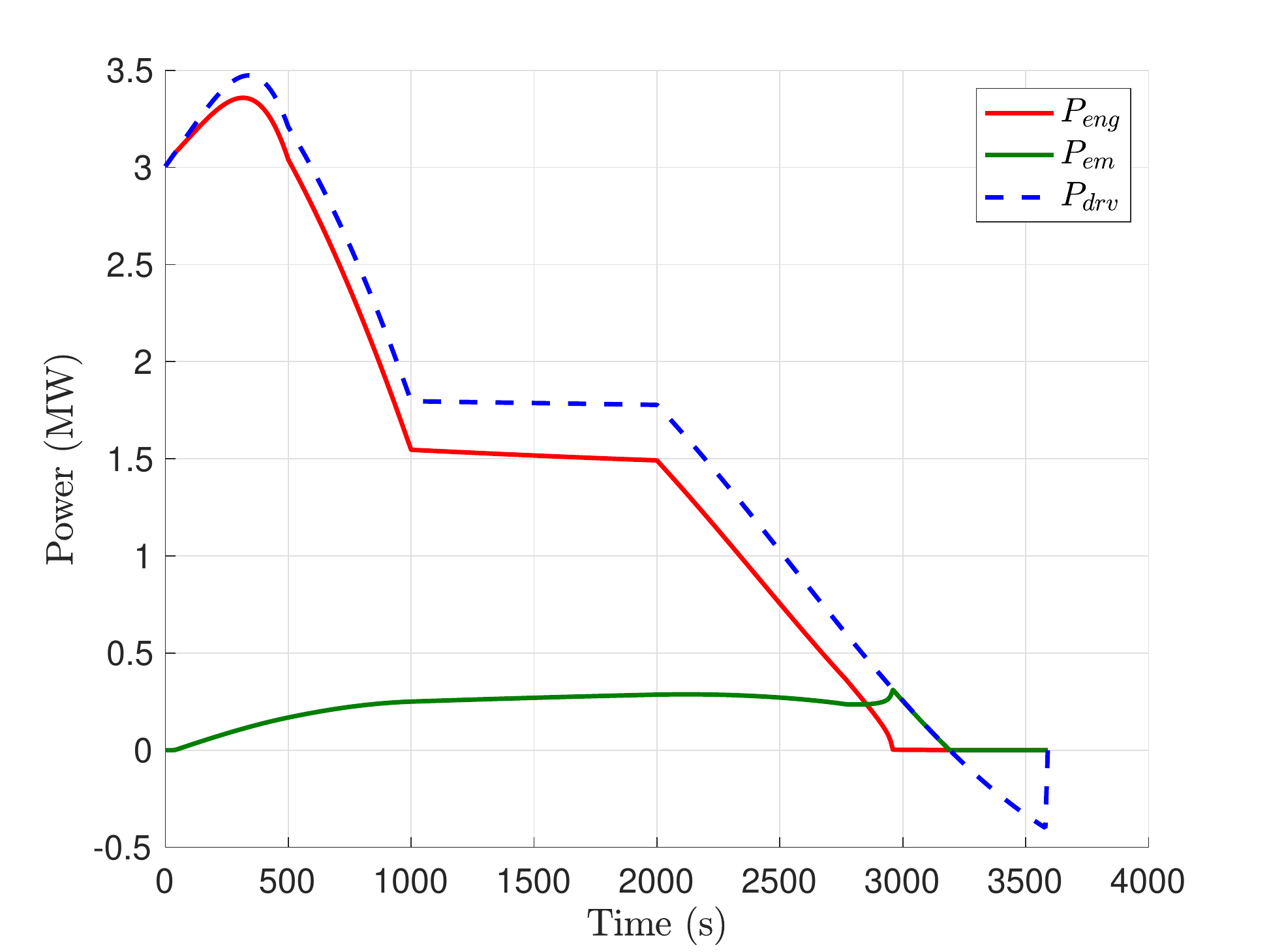}
\caption{Optimal power split for the case of a fuel map with an increased rate of fuel consumption.}
\label{fig:Psplit2}
\end{figure}

More radical optimal power split solutions are obtained if the change in aircraft mass during flight is more significant. In particular, the MPC strategy allocates more electrical power at the end of the flight if the gas turbine fuel consumption is increased. For example, Figure \ref{fig:Psplit2} shows the power split solution for a situation in which the rate of fuel consumption is increased so that the change in aircraft total mass during flight is 15\% (with all other simulation parameters unchanged). 
%
%
The fuel consumption for the CDCS strategy in this case is about $4\%$ higher than that of the MPC strategy. 

\begin{figure}
\includegraphics[width = .45\textwidth]{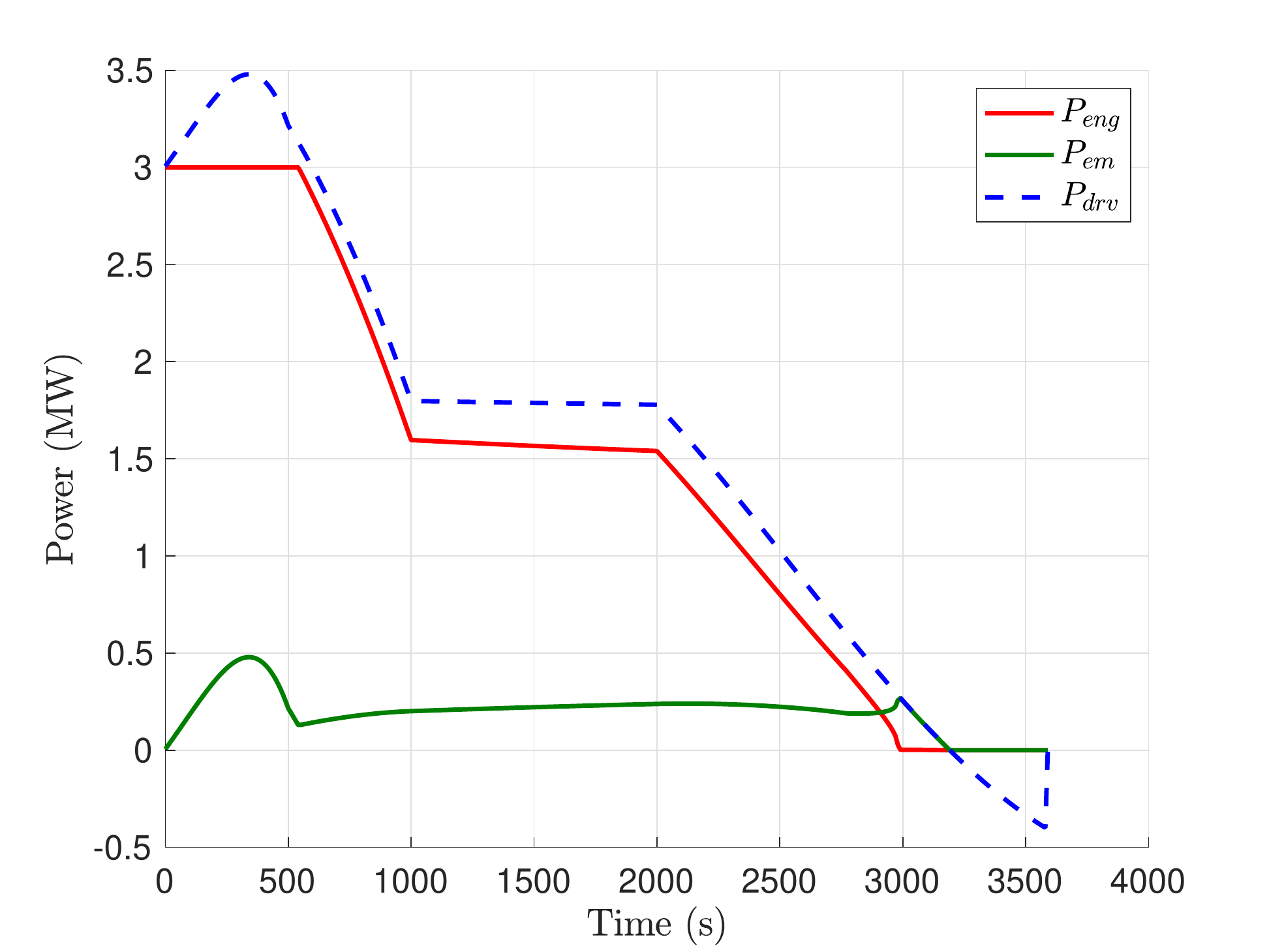}
\caption{Optimal power split with gas turbine saturation.}
\label{fig:Psplit3}
\end{figure}


We next consider the case that the upper limit on the gas turbine power output is reduced to $\overline{P}_{\gt}=3$ $\si{MW}$. The MPC energy management strategy for this case is shown in Figure \ref{fig:Psplit3}. 
Here the power demand is such that the gas turbine is at maximum power while the aircraft climbs. As a result, the electric motor is needed to meet the total power output requirement while the gas turbine power output is saturated.
The fuel consumption for this scenario is increased slightly (by $0.1\%$) since the electrical power is mostly used at the beginning of the flight to compensate for the limit on the gas turbine power output. 
%


\section{Conclusions}
\label{sec:conclusion}

This paper proposes a model predictive control law for energy management in hybrid-electric aircraft.
The main contribution of the work is a convex formulation of the problem of minimising fuel consumption for a given future flight path.
We provide a simulation study to illustrate the approach, and demonstrate that significant fuel savings can be achieved relative to heuristic strategies.
The convexity of the formulation is crucial for computational tractability and is expected to be a basic requirement for verification by the aviation industry.
Future work will consider the design of bespoke solvers. In particular, first order solution methods are expected to provide computational savings by exploiting the high degree of separability in the problem, while also being suitable for real-time implementation.
%
The modelling approach described in this paper provides a framework for optimising system design, and future work will explore flight path optimisation and evaluate alternative hybrid propulsion configurations.
%


\bibliography{convex_hybrid_aircraft_ifac}

\begin{thebibliography}{14}
\providecommand{\natexlab}[1]{#1}
\providecommand{\url}[1]{\texttt{#1}}
\providecommand{\urlprefix}{URL }
\expandafter\ifx\csname urlstyle\endcsname\relax
  \providecommand{\doi}[1]{doi:\discretionary{}{}{}#1}\else
  \providecommand{\doi}{doi:\discretionary{}{}{}\begingroup
  \urlstyle{rm}\Url}\fi

\bibitem[{Abbott et~al.(1945)Abbott, Von~Doenhoff, and Stivers~Jr}]{abbott1945}
Abbott, I., Von~Doenhoff, A., and Stivers~Jr, L. (1945).
\newblock Summary of airfoil data.
\newblock Technical report, Langley Memorial Aeronautical Laboratory.

\bibitem[{Di~Cairano et~al.(2014)Di~Cairano, Bernardini, Bemporad, and
  Kolmanovsky}]{dicairano14}
Di~Cairano, S., Bernardini, D., Bemporad, A., and Kolmanovsky, I. (2014).
\newblock Stochastic {MPC} with learning for driver-predictive vehicle control
  and its application to {HEV} energy management.
\newblock \emph{IEEE Trans.\ Control Syst.\ Technol.}, 22(3), 1018--1031.

\bibitem[{East and Cannon(2018)}]{east2018}
East, S. and Cannon, M. (2018).
\newblock An {ADMM} algorithm for {MPC}-based energy management in hybrid
  electric vehicles with nonlinear losses.
\newblock In \emph{2018 IEEE Conference on Decision and Control (CDC)},
  2641--2646.

\bibitem[{East and Cannon(2019)}]{east19ieeetcst}
East, S. and Cannon, M. (2019).
\newblock Energy management in plug-in hybrid electric vehicles: Convex
  optimization algorithms for model predictive control.
\newblock \emph{IEEE Trans.\ Control Syst.\ Technol. (Early Access)}.

\bibitem[{Grant and Boyd(2008)}]{cvx}
Grant, M. and Boyd, S. (2008).
\newblock Graph implementations for nonsmooth convex programs.
\newblock In V.~Blondel, S.~Boyd, and H.~Kimura (eds.), \emph{Recent Advances
  in Learning and Control}, 95--110. Springer-Verlag Limited.

\bibitem[{Hall et~al.(2017)Hall, Huang, Uranga, Greitzer, Drela, and
  Sato}]{hall2017}
Hall, D.K., Huang, A.C., Uranga, A., Greitzer, E.M., Drela, M., and Sato, S.
  (2017).
\newblock Boundary layer ingestion propulsion benefit for transport aircraft.
\newblock \emph{Journal of Propulsion and Power}, 33(5), 1118--1129.

\bibitem[{Josevski et~al.(2017)Josevski, Katriniok, and Abel}]{josevski17}
Josevski, M., Katriniok, A., and Abel, D. (2017).
\newblock Scenario {MPC} for fuel economy optimization of hybrid electric
  powertrains on real-world driving cycles.
\newblock In \emph{Proc.\ American Control Conference}, 5629--5635.

\bibitem[{Kim et~al.(2011)Kim, Cha, and Peng}]{kim11}
Kim, N., Cha, S., and Peng, H. (2011).
\newblock Optimal control of hybrid electric vehicles based on {P}ontryagin's
  {M}inimum {P}rinciple.
\newblock \emph{IEEE Trans.\ Control Syst.\ Technol.}, 19(5), 1279--1287.

\bibitem[{Koot et~al.(2005)Koot, Kessels, de~Jager, Heemels, van~den Bosch, and
  Steinbuch}]{koot05}
Koot, M., Kessels, J., de~Jager, B., Heemels, W., van~den Bosch, P., and
  Steinbuch, M. (2005).
\newblock Energy management strategies for vehicular electric power systems.
\newblock \emph{IEEE Trans.\ Veh.\ Technol.}, 54(3), 771--782.

\bibitem[{Lin et~al.(2003)Lin, Peng, Grizzle, and Kang}]{lin03}
Lin, C., Peng, H., Grizzle, J., and Kang, J. (2003).
\newblock Power management strategy for a parallel hybrid electric truck.
\newblock \emph{IEEE Trans.\ Control Syst.\ Technol.}, 11(6), 839--849.

\bibitem[{Onori and Tribioli(2015)}]{onori15}
Onori, S. and Tribioli, L. (2015).
\newblock Adaptive {P}ontryagin's {M}inimum {P}rinciple supervisory controller
  design for the plug-in hybrid {GM} {C}hevrolet {V}olt.
\newblock \emph{Appl.\ Energy}, 147, 224--234.

\bibitem[{Sciarretta and Guzzella(2007)}]{sciarretta07}
Sciarretta, A. and Guzzella, L. (2007).
\newblock Control of hybrid electric vehicles.
\newblock \emph{IEEE Control Syst.\ Mag.}, 27(2), 60--70.

\bibitem[{Stevens et~al.(2016)Stevens, Lewis, and
  Johnson}]{stevens2015aircraft}
Stevens, B.L., Lewis, F.L., and Johnson, E.N. (2016).
\newblock \emph{Aircraft Control and Simulation: Dynamics, Controls Design, and
  Autonomous Systems}.
\newblock John Wiley \& Sons.

\bibitem[{Thipphavong et~al.(2018)}]{nasa_urban}
Thipphavong, D.P. et~al. (2018).
\newblock Urban air mobility airspace integration concepts and considerations.
\newblock In \emph{AIAA Aviation Forum (Aviation 2018)}. Atlanta, GA.

\end{thebibliography}

\end{document}